\documentclass[12pt]{article}

\usepackage{amsfonts}
\usepackage{xr}
\usepackage{graphicx}
\usepackage{amsmath}

\def\to{\rightarrow}

\def\R{\mathbb{R}}

\def\tilde{\widetilde}
\def\epsilon{\varepsilon}
\def\trait (#1) (#2) (#3){\vrule width #1pt height #2pt depth #3pt}
\def\fin{\hfill\trait (0.1) (5) (0) \trait (5) (0.1) (0) \kern-5pt \trait (5) (5) (-4.9) \trait (0.1) (5) (0)}

\newenvironment{note}[1]{\footnote{\footnotesize{#1}}}

\newtheorem{thm}{\bf Theorem}[section]
\newtheorem{lem}[thm]{\bf Lemma}
\newtheorem{prop}[thm]{\bf Proposition}

\title{Some dependence results between the spreading speed and the coefficients of the space-time {Fisher-KPP} equation}
\date{}
\author{Gregoire Nadin\thanks{
D\'epartement de Math\'ematiques et Applications,
\'Ecole Normale Sup\'erieure, CNRS UMR8553 ,
            45 rue d'Ulm, F~75230 Paris cedex 05, France ; e-mail: nadin@dma.ens.fr}}

\begin{document}
\maketitle

\begin{abstract}
 We investigate in this paper the dependence relation between the space-time periodic coefficients $A,q$ and $\mu$ of the reaction-diffusion equation
$$\partial_t u -\nabla \cdot (A(t,x)\nabla u)+q(t,x)\cdot \nabla u = \mu(t,x)u(1-u),$$
and the spreading speed of the solutions of the Cauchy problem associated with this equation and compactly supported initial data. 
We prove in particular that
(1) taking the spatial or temporal average of $\mu$ decreases the minimal speed, 
(2) if the coefficients do not depend on $t$ and $q\equiv 0$, then increasiong the amplitude of the diffusion matrix $A$ increases the minimal speed,
(3) if $A=I_N$, $\mu$ is a constant, then the introduction of a space periodic drift term $q=\nabla Q$ increases the minimal speed.

To prove these results, we use a variational characterization of the spreading speed that involves a family of periodic principal eigenvalues associated with the linearization of 
the equation near $0$. We are thus back to the investigation of the dependence relation between this family of eigenvalues and the coefficients. 
\end{abstract}


\bigskip


\tableofcontents


\section{Introduction}


\subsection{General framework and definition of the spreading speed}

This article investigates the asymptotic properties of the solutions of the space-time periodic Fisher-KPP equation:
\begin{equation} \label{eqprinc}\left\{ \begin{array}{l}
\partial_t u -\nabla \cdot (A(t,x)\nabla u)+q(t,x)\cdot \nabla u = \mu(t,x)u(1-u) \hbox{ in } \R^+\times \R^N,\\
u(0,x)=u_0(x) \hbox{ in } \R^N,\\ \end{array}
\right.\end{equation}
where $u_0$ is a nonnegative, continuous and compactly supported initial datum. This equation arises in various models, that comes from genetics, population dynamics, 
combustion, chemistry etc. In these models, the function $u$ represents a density of population or of a chemical material. It diffuses in a space-time heterogeneous media 
through a diffusion matrix $A(t,x)$ and it reacts through a reaction term $\mu(t,x)u(1-u)$, where $\mu(t,x)$ represents a growth rate at small density. Lastly, 
it is advected at a speed $q(t,x)$. 

This equation has first been investigated in one-dimensional media by Kolmogorov, Petrovski and Piskunov \cite{KPP} and by Fisher \cite{Fisher} in the 30's, 
then in multidimensional media by Aronson and Weinberger \cite{Aronson} in the 70's, under the hypotheses $A=a I_N$ and $\mu$ do not depend on $(t,x)$, $a>0$, $\mu>0$ and $q\equiv 0$
Among other properties, these authors proved that for all $e\in\mathbb{S}^{N-1}$, if $u$ is the solution of (\ref{eqprinc}) and $u_0$ is compactly supported,
\begin{equation} \label{spreading-hom} \left\{\begin{array}{rl}
\displaystyle{\mathop{\liminf}_{t\to+\infty}}\ u(t,cte)=1 & \hbox{ if }0\le c<2\sqrt{\mu a},\vspace{4pt}\\
\displaystyle{\mathop{\lim}_{t\to+\infty}} u(t,cte)=0 & \hbox{ if }c>2\sqrt{\mu a}.\end{array}\right.
\end{equation}
This result is called a {\em spreading property} and the speed $c^*=2\sqrt{\mu a}$ is called a {\em spreading speed} in direction $e$. 

In the late 70's, spreading properties have been proved for the space periodic Fisher-KPP equation by Freidlin and Gartner \cite{Gartner} and Freidlin \cite{Freidlin2}. In such
media the spreading speed depend on the direction of propagation $e$ and is defined in an implicit way. Namely, these authors proved that, 
if $\mu>0$, $\nabla\cdot q=0$, $q$ has a null average and $A$, $q$ and $\mu$ do 
not depend on $t$, then for all $e\in\mathbb{S}^{N-1}$, there
exists a speed $c^*_e$ such that 
\begin{equation} \label{spreading} \left\{\begin{array}{rl}
\displaystyle{\mathop{\liminf}_{t\to+\infty}}\ u(t,cte)=1 & \hbox{ if }0\le c<c^*_e,\vspace{4pt}\\
\displaystyle{\mathop{\lim}_{t\to+\infty}} u(t,cte)=0 & \hbox{ if }c>c^*_e.\end{array}\right.
\end{equation}
Moreover, they proved a useful variational formula for $c^*_e$ that we will give in Section \ref{section-eigen}. 
The existence of pulsating traveling front\footnote{We refer to the references below for the definition of this notion.} for the space periodic Fisher-KPP equation has later been proved under various hypotheses by Xin \cite{Xin}, Berestycki and Hamel \cite{Frontexcitable}
and Berestycki, Hamel and Roques \cite{Base2}. These authors found a link between the spreading speed in direction $e$ and the minimal speed of existence of pulsating traveling fronts
in the directions $\xi$ such that $e\cdot \xi>0$. 

The investigation of the dependence relations between the coefficients $(A,q,\mu)$ of the space periodic Fisher-KPP equation and the spreading speed 
$c^*_e=c^*_e(A,q,\mu)$ has started at the beginning of the 
2000's. There is a wide litterature on this topic, that we will describe in Section \ref{section-results} below. 

Recently, spreading properties in space-time periodic media have been proved when $N=1$, $A=I_N$, $\mu>0$ is a constant, $\nabla\cdot q$ and $q$ has a null average
by Nolen, Rudd and Xin \cite{NolenRuddXin} and 
under the general hypotheses of Section \ref{section-hyp} below by Weinberger \cite{Weinberger} and by Berestycki,
Hamel and Nadin \cite{BHNa}. Hardly no dependence results between the spreading speed and the coefficients have been obtained in this case. 

\bigskip

In the present paper we give new dependence results between the spreading speed $c^*_e=c^*_e(A,q,\mu)$ and the space-time coefficients $(A,q,\mu)$. Some
of these results are extension of results that are known from space periodic media to space-time periodic media. But we also prove some dependence relations
that are new, even in space periodic media. Our main results are the following:
\begin{itemize}
\item taking the spatial or temporal average of $\mu$ decreases the minimal speed, 
 \item if the coefficients do not depend on $t$ and $q\equiv 0$, then increasing the amplitude of the diffusion matrix $A$ increases the minimal speed,
 \item if $A=I_N$, $\mu$ is a constant, then the introduction of a space periodic drift term $q=\nabla Q$ increases the minimal speed.
\end{itemize}

\smallskip

{\bf Organization of the paper.} In Section \ref{section-hyp}, we give the hypotheses we require on the coefficients $(A,q,\mu)$. Then, we define the family of periodic principal eigenvalues
involved in the variational characterization
of the spreading speed and we clearly state this characterization in Section \ref{section-eigen}. We state our results in Section \ref{section-results}. We also give a review of the known 
dependence relations between the coefficients $(A,q,\mu)$ and the spreading speed in this Section. Lastly, we prove our dependence relations with respect to $\mu$
in Section \ref{section-proofsmu}, with respect to $A$ in Section \ref{section-proofsA} and with respect to $q$ in Section \ref{section-proofsq}.


\subsection{Hypotheses} \label{section-hyp}

We assume that the diffusion matrix $A$, the advection term $q$ and the growth rate $\mu$ are periodic in $(t,x)$. 
That is, there exist some positive constant $T$ and some vectors $L_{1},...,L_{N}$, where $L_{i}$ is colinear to the axis of coordinates $e_{i}$, such that for all 
$i\in [1,N]$, for all $(t,x)\in\R\times\R^{N}$, one has:
\[\begin{array}{ccccc}
A(t,x+L_{i})&=&A(t+T,x)&=&A(t,x),\\
\mu(t,x+L_{i})&=&\mu(t+T,x)&=&\mu(t,x),\\
q(t,x+L_{i})&=&q(t+T,x)&=&q(t,x).\\
\end{array}\]
We define the periodicity cell $C= \Pi_{i=1}^{N} (0,|L_{i}|)$. In the sequel the notion of periodicity will always refer to the periods $(T,L_1,...,L_N)$. 

We shall need some regularity assumptions on $\mu, A, q$.
The growth rate $\mu: \R\times\R^N\to \R$ is supposed to be of class $C^{\frac{\delta}{2},\delta}$
The matrix field $A:\mathbb{R} \times \mathbb{R}^{N}
\rightarrow S_{N}(\mathbb{R})$ is supposed to be of class $C^{\frac{\delta}{2},1+\delta}$. We suppose futhermore that $A$ is uniformly elliptic and continuous: there exist some positive constants $\gamma$ and $\Gamma$ such that for all $\xi \in \mathbb{R}^{N}, (t,x) \in \mathbb{R} \times \mathbb{R}^{N}$ one has:
\begin{equation} \label{ellipticity-dep} \gamma |\xi|^{2} \leq \sum_{1 \leq i,j \leq N} a_{i,j}(t,x) \xi_{i} \xi_{j} \leq \Gamma |\xi|^{2}, \end{equation}
where $|\xi|^{2}=\xi_{1}^{2}+...+\xi_{N}^{2}$. and $a_{i,j}(t,x)$ is the coefficient $(i,j)$ of the matrix $A(t,x)$. 
The drift term $q: \mathbb{R} \times \mathbb{R}^{N} \rightarrow \mathbb{R}^{N}$ is supposed to be of class $C^{\frac{\delta}{2},\delta}$ and we assume that $\nabla\cdot q\in L^{\infty}(\R\times\R^{N})$.
In the sequel, the direction of propagation $e\in\mathbb{S}^{N-1}$ will be fixed. 


\subsection{Characterization of the spreading speed with periodic principal eigenvalues} \label{section-eigen}

The characterization of the spreading speed invloves the family of operators which is associated with exponentially decreasing solutions of the linearization of 
(\ref{eqprinc}) in the neighborhood of $0$:
\begin{equation} L_{\lambda} \psi = \partial_{t}\psi - \nabla\cdot (A\nabla \psi) -2\lambda A \nabla \psi+q\cdot\nabla\psi-(\lambda A\lambda+\nabla\cdot(A\lambda)+
\mu-q\cdot\lambda)\psi,\end{equation}
where $\lambda \in \mathbb{R}^{N}$ and $\psi \in \mathcal{C}^{1,2}(\mathbb{R} \times \mathbb{R}^{N})$. 
Then it has been proved by the author in \cite{eigenvaluearticle} that there exists a unique real number $k_{\lambda}(A,q,\mu)$ such that 
there exists a function $\psi \in \mathcal{C}^{1,2}(\mathbb{R} \times \mathbb{R}^{N}) $ that satisfies
\begin{equation} \label{eigentw2} \left\{\begin{array}{l}
L_{\lambda} \psi = k_\lambda(A,q,\mu) \psi \ \hbox{in} \ \R\times\R^N,\\
\psi > 0 \ \hbox{in} \ \R\times\R^N,\\
\psi \ \hbox{is periodic}.\\
\end{array}\right.
\end{equation}
We call $k_\lambda(A,q,\mu)$ the {\em space-time periodic principal eigenvalue} associated with operator $L_\lambda$. 

The variational characterization of the spreading speed we will use in the sequel has been proved by Berestycki, Hamel and the author \cite{BHNa} in two different ways:
\begin{thm}\cite{BHNa}
Assume that $k_\lambda(A,q,\mu)<0$ for all $\lambda\in\R^N$. Then if $u_0\not\equiv 0$ is a compactly supported, continuous and nonnegative initial datum and $u$ is the 
associated solution of the Cauchy problem (\ref{eqprinc}), one has
 \begin{equation}  \left\{\begin{array}{rl}
\displaystyle{\mathop{\liminf}_{t\to+\infty}}\ u(t,cte)=1 & \hbox{ if }0\le c<c^*_e(A,q,\mu),\vspace{4pt}\\
\displaystyle{\mathop{\lim}_{t\to+\infty}} u(t,cte)=0 & \hbox{ if }c>c^*_e(A,q,\mu),\end{array}\right.
\end{equation}
with 
\begin{equation} \label{charac-c} c^*_e(A,q,\mu)=\min_{\lambda\cdot e<0} \frac{k_{\lambda}(A,q,\mu)}{\lambda\cdot e}.\end{equation}
\end{thm}

This formula highly simplifies the investigation of the dependence relation between $(A,q,\mu)$ and $c^*_e(A,q,\mu)$.


\section{Statement of the dependence results} \label{section-results}


\subsection{Spatial and temporal averaging of the growth rate}

We begin with two comparison principles with the averaged media in $x$ or in $t$:

\begin{prop} \label{dep-x} {\bf (Influence of the spatial variations)}
If $A$ and $q$ do not depend on $x$, define
\[\overline{\mu}(t)=\frac{1}{|C|}\int_{C}\mu(t,x)dx.\] 
Then, if $\int_{(0,T)\times C}\mu\geq 0$, the following comparison holds: 
\begin{equation} \label{dep-x-ineq} c^{*}_e(A,q,\mu)\geq c^{*}_e(A,q,\overline{\mu})=_min_{e\cdot \xi>0}\frac{2}{T}\sqrt{\int_{0}^{T}\xi A \xi\int_{0}^{T}\mu}-\frac{1}{T}\int_{0}^{T}q\cdot \xi.\end{equation}
Moreover, the equality holds if and only if $\mu$ does not depend on $x$.
\end{prop}

This theorem means that, somehow, the heterogeneity in $x$ of the growth rate increases the speed of propagation. 
Using this heuristic definition of "heterogeneity", it is not true that heterogeneous drift or diffusion coefficient speed up the propagation. 
We will prove later that some compressible drifts may slow down the propagation. It has also been proved by Papanicolaou and Xin \cite{PapanicolaouXin} that, 
in dimension $1$, if $b$ is a space periodic continuous function of average $0$ and $\mu_0$ is a positive constant, then 
$c^*_{e_1}(1+\delta b,0,\mu_0)\leq c^*_{e_1}(1,0,\mu_0)=2\sqrt{\mu_0}$ when $\delta$ is small enough. 
Taking the additive average of the diffusion or advection coefficients is not the good mean to quantify the heterogeneity. 
Other kind of averaging may give positive result. For example, it has been proved by the author \cite{Steiner} that, for space periodic media, 
in dimension $1$, $c^*_{e_1}(<a>_H,0,\overline{\mu})\leq c^*_{e_1}(a,0,\mu)$, where $<a>_H$ is the harmonic average of $a$. 

\bigskip

Similarly, taking the temporal average of the growth rate decreases the minimal speed of propagation:
\begin{thm} \label{influencetimevar}{\bf (Influence of the temporal variations)}
If $A$ and $q$ do not depend on $t$, define
\[\hat{\mu}(x)=\frac{1}{T}\int_0^T\mu(t,x)dt.\] 
Then, if $k_0(A,q,\hat{\mu})<0$, the following comparison holds:
\[c^{*}_e(A,q,\mu)\geq c^{*}_e(A,q,\hat{\mu}).\]
Moreover, the equality holds if and only if $\mu$ can be written:
$\mu(t,x)=\mu_{1}(x)+\mu_{2}(t).$
\end{thm}

\subsection{Influence of the amplitude of the reaction term}

We first state that increasing the reaction term increases the speed of propagation. This is an easy extension of Proposition 1.15 of \cite{SpeedKPP}. 

\begin{prop} \label{propamp}
If $\mu_{1}\geq \mu_{2}$, then for all $A,q$, one has:
\[c^*_e(A,q,\mu_{1})\geq c^*_e(A,q,\mu_{2}).\]
Moreover, the equality holds if and only if $\mu_{1}\equiv \mu_{2}$.
\end{prop}

Next, one can wonder what is the influence of the amplitude of the growth rate on the minimal speed. 

\begin{prop} \label{dep-amp-classical}

1. Assume that $\mu_0$ is a constant and consider a space-time periodic function $\eta$.
If $\int_{(0,T)\times C} \eta\geq 0$ (resp. $\int_{(0,T)\times C} \eta> 0$), then $B\mapsto
c^*_e(I_N,0,\mu_0+B\eta)$ is nondecreasing (resp. increasing). Moreover, if $\int_{(0,T)\times C} \eta\geq 0$, then $B\mapsto c^*_e(I_N,0,\mu_0+B\eta)$ is increasing if and only if $\eta$ is not a constant with respect to $x$.

2. Assume that $A$, $q$ and $\mu$ do not depend on $t$ and that $\max_{x\in\R^{N}}\int_0^T\eta(t,x)dt >0$. Then $B\mapsto
c^*_e(A,q,\mu+B\eta)$ is increasing for $B$ large enough.
\end{prop}

This result extends that of Berestycki, Hamel and Roques \cite{Base2} from space periodic to space-time periodic media. 
In such media, the hypothesis of 2. only reads $\max_{\R^N}\eta>0$. Hence, the main interest of Proposition \ref{dep-amp-classical} is to identify the generalization of this hypothesis 
to time-dependent media, that is, $\max_{x\in\R^{N}}\int_0^T\eta(t,x)dt >0$.

\subsection{Montonicity with respect to the diffusion amplitude}

It seems natural that increasing the diffusion coefficient may increase the speed of propagation. In \cite{SpeedKPP}, 
Berestycki, Hamel and Nadirashvili have proved that $\kappa\mapsto c^{*}_e(\kappa A,0,\mu)$ is increasing if $\mu$ is constant 
and $A$ only depends on $x$. It was an open problem to generalize this result to heterogeneous growth rate. 

\begin{thm} \label{mon-diffusion}{\bf (Monotonicity with respect to the diffusion in space periodic media)} 
 Assume that $A$ and $\mu$ do not depend on $t$ and that $q\equiv 0$, then 
$$\kappa\mapsto c^*_e(\kappa A, 0, \mu) \hbox{ is increasing}.$$
\end{thm}

Unfortunately, such a generalization is not always true. If one includes a drift term that depends on $x$, El Smaily has proved in \cite{ElSmaily}, that the associated speed is not monotonic in $\kappa$ in general. El Smaily has also
proved that $A\geq B$ in the sense of positive matrix does not imply $c^*_e(A,0,\mu)\geq c^*_e(B,0,\mu)$.  

In order to conclude this section, let us mention some dependence results for the function $\kappa\mapsto \displaystyle\frac{c^*_e(\kappa A, 0, \mu)}{\sqrt{\kappa}}$.
The author has proved that this function is nonincreasing in \cite{Steiner}. He has also computed the limit of 
this function when $\kappa\rightarrow +\infty$. In dimension $1$, the limit when $\kappa\to 0$ has been computed by Hamel, Roques and Fayard \cite{HamelRoqueslimgamma} when $A$ and $\mu$ only take two values
and by Hamel, Roques and the author \cite{limgamma} for general $A$ and $\mu$ with $q\not\equiv 0$.

\subsection{Influence of the drift}

\subsubsection*{Incompressible drifts}

It has been proved by Berestycki, Hamel and Nadirashvili \cite{SpeedKPP} that, in space periodic media, the introduction of an incompressible drift with null average 
increases the propagation speed. Actually, the difficulty is to understand what it is the amplitude of this speed-up. 
It has been proved that this speed-up depends on the geometric properties of the level-lines of the flow associated with $q$ 
(see \cite{Audoly, Largedrift, Heinze2005, Kiselevadvection, RyzhikZlatos, Zlatosc/A, Zlatoscex}).

We only consider here the case of a shear flow. 

\begin{prop} \label{shearflowprop}
Assume that $\mu_0$ is a positive constant and that the drift term can be written $q(t,x)=(q_{1}(t,y),0,...,0)$, where one writes 
$x=(x_{1},y)\in\R\times\R^{N-1}$, $q_1\not\equiv 0$ and $\int_C q_1=0$. 
Then for all $e\in\mathbb{S}^{N-1}$, $B\mapsto c^*_e(I_N,Bq,\mu_0)$ is increasing.
\end{prop}

This monotonicity has been numerically observed by Nolen and Xin \cite{Nolenshear} in the case $e=e_1$. There was no analytical proof of this numerical observation before, 
as far as we know. 

\subsubsection*{Compressible drifts}

There is hardly no paper on the influence of a compressible drift on the speed in the litterature. 
Only Nolen and Xin have investigated propagation problems in such media before \cite{NolenXin1d}. 
For space stationary random drifts, when $N=1$, $A=1$ and $\mu_0$ is a positive constant. 
They have proved that:
\[ \forall \alpha \in (0,1), \ \exists c, \ C_\alpha, \ \forall B>0, \frac{c}{B}\leq c^*_{e_1}(1,Bq,\mu_0)\leq \frac{C_\alpha}{B^\alpha},\]
where $c^*_e(1,Bq,\mu_0)$ is the spreading speed in direction $e_1$ associated with the drift $Bq$. 
We focus here on the drifts that can be written $q=A\nabla Q$ and prove that such drifts slow down the propagation. 

\begin{thm} \label{potdrift} {\bf (Influence of a drift $q=A\nabla Q$)}
Assume that $A=I_N$, $\mu_0$ is a positive constant, $q$ does not depend on $t$, $\int_C q =0$ and that $q$ can be written $q=\nabla Q$. Then:

1. $\displaystyle \frac{c^*_e(I_N,B\nabla Q,\mu_0)}{B}\rightarrow 0 \ \hbox{as} \ B\rightarrow+\infty,$

2. one has $$c^*_e(I_N,B\nabla Q,\mu_0)\leq c^*_e(I_N,0,\mu_0)=2\sqrt{\mu_0}.$$
\end{thm}

\noindent {\bf Remark.} In dimension $1$, if $A=I_{N}$, the hypothesis is equivalent to $\int_{(0,T)\times C}q=0$. In dimension $2$ or $3$, if $A=I_{N}$, 
it is equivalent to $\int_{(0,T)\times C}q=0$ and $curl \ q=0$.


\section{Proof of the dependence results with respect to the growth rate} \label{section-proofsmu}

The aim of this section is to prove Proposition \ref{dep-x}, Theorem \ref{influencetimevar}, Proposition \ref{propamp} and Proposition \ref{dep-amp-classical}. 
We first give a direct proof of Proposition \ref{propamp}. Then we state some general dependence relations with respect to the growth rate that will enable us to prove 
the other results. 

\subsection{Proof of Proposition \ref{propamp}}

\noindent {\bf Proof of Proposition \ref{propamp}.} We use the same kind of proof as Berestycki, Hamel and Nadirashvili in \cite{SpeedKPP}. 
If $\mu_{1}\geq \mu_{2}$, one immediatly gets $k_{\lambda}(A,q,\mu_{1})\leq k_{\lambda}(A,q,\mu_{2})$ using the min-max characterization of $k_{\lambda}(A,q,\mu)$ proved
by the author in \cite{eigenvaluearticle}:
\begin{equation} \label{minmax-dep} k_{\lambda}(A,q,\mu)= \max_{\phi\in\mathcal{C}^{1,2}_{per}(\R\times\R^N), \ \phi>0.}\min_{\mathbb{R} \times \mathbb{R}^{N}} \Big( \frac{L_\lambda\phi}{\phi}\Big).\end{equation}
Thus $c^*_e(A,q,\mu_{1})\geq c^*_e(A,q,\mu_{2})$. Assume now that $c^*_e(A,q,\mu_{1})= c^*_e(A,q,\mu_{2})$ and take some $\lambda\in\R^N$ such that $\lambda\cdot e<0$ and
\[c^*_e(A,q,\mu_{1})=\displaystyle \frac{k_{\lambda}(A,q,\mu_{1})}{\lambda\cdot e}.\] 
One has:
\[c^*_e(A,q,\mu_{1})=\frac{k_{\lambda}(A,q,\mu_{1})}{\lambda\cdot e}=c^*_e(A,q,\mu_{2})\leq \frac{k_{\lambda}(A,q,\mu_{2})}{\lambda\cdot e},\]
and then $k_{\lambda}(A,q,\mu_{1})\geq k_{\lambda}(A,q,\mu_{2})$, that is, $k_{\lambda e}(A,q,\mu_{1})= k_{\lambda e}(A,q,\mu_{2})$. 

Take now $\phi_{\lambda}^{1}$ some eigenfunction associated with $\mu_{1}$ and $\phi_{\lambda}^{2}$ some eigenfunction associated with $\mu_{2}$. 
Set $\kappa=\max_{(0,T)\times C}\displaystyle\frac{\phi_{\lambda}^{1}}{\phi_{\lambda}^{2}}$ and $z=\phi_{\lambda}^{1}-\kappa\phi_{\lambda}^{2}$. This function is nonpositive, vanishes somewhere and satisfies:
\[\begin{array}{l}\partial_{t}z - \nabla\cdot (A\nabla z) -2\lambda eA \nabla z+q\cdot\nabla z-(\lambda^{2} eAe+\lambda\nabla\cdot(Ae)+\mu_{1}-\lambda q\cdot e+k_{\lambda e}(A,q,\mu_{1}))z\\
=(\mu_{2}-\mu_{1})\phi_{\lambda}^{2}\leq 0.\\\end{array}\]
Thus the periodicity in $t$ and the strong parabolic maximum principle give $z\equiv 0$. Hence $\mu_{1}\equiv \mu_{2}$.$\Box$


\subsection{Strict concavity of the principal eigenvalue}

The concavity of $\mu\mapsto k_\lambda(A,q,\mu)$ has already been proved in \cite{eigenvaluearticle}, but we focus here on the strict concavity, which will be our main tool in order to  investigate equality cases later.

\begin{prop} \label{concavitydep}
 For all $A$, $q$, $\mu_1$, $\mu_2$, $r\in (0,1)$ and $\lambda\in\R^N$, one has
\[k_\lambda(A,q,r\mu_1+(1-r)\mu_2)\geq r k_\lambda(A,q,\mu_1)+(1-r) k_\lambda(A,q,\mu_2).\]
Moreover, if $A$ and $q$ do not depend on $t$, the equality holds if and only if $\mu_1-\mu_2$ does not depend on $x$.
\end{prop}

\noindent {\bf Proof.} 
As we already mentionned it, the concavity has already been proved by the author in \cite{eigenvaluearticle}. We include this proof here by sake of completeness and 
because it will lead us to the strict concavity. 

As we are considering any possible $q$, $\mu_1$ and $\mu_2$, we can assume\note{This is where the hypothesis ``$A$ and $q$ do not depend on $x$'' is used in the equality case.} that $\lambda=0$.
Set $\mu = r\mu_{1} + (1-r)\mu_{2}$ and consider $\phi_{1}$ and $\phi_{2}$ some periodic principal eigenfunctions associated with $\mu_1$ and $\mu_2$. Define $z_{1} = \ln(\phi_{1})$,
$z_{2} = \ln(\phi_{2})$, $z = rz_{1} + (1-r) z_{2}$ and $\phi = e^{z}$. One can compute:
\[ \frac{\partial_{t} \phi -\nabla\cdot(A\nabla\phi)+q\cdot\nabla\phi}{\phi} = \partial_{t}z-\nabla\cdot(A\nabla z) - \nabla z A \nabla z+q\cdot \nabla z,\]
\begin{equation} \label{strictconc0} \begin{array}{rcl} \hbox{and } \nabla z A \nabla z &=& r \nabla z_{1} A \nabla z_{1} + (1-r) \nabla z_{2} A \nabla z_{2} -r(1-r)(\nabla z_{1} - \nabla z_{2})A(\nabla z_{1} - \nabla z_{2})\\
&\leq& r \nabla z_{1} A \nabla z_{1} + (1-r) \nabla z_{2} A \nabla z_{2}.  \\
\end{array} \end{equation}

Hence, for all $(t,x)\in\R\times\R^N$:
\[ \begin{array}{rcl}\displaystyle\frac{\partial_{t} \phi -\nabla\cdot(A\nabla\phi)+q\cdot\nabla\phi}{\phi} - \mu &\geq& r(\partial_{t}z_{1}-\nabla\cdot(A\nabla z_{1}) - \nabla z_{1} A \nabla z_{1} +q\cdot \nabla z_{1}- \mu_{1})\\
&&+(1-r)(\partial_{t}z_{2}-\nabla\cdot(A\nabla z_{2}) - \nabla z_{2} A \nabla z_{2}+q\cdot\nabla z_{2}-\mu_{2}) \\
&&\\
&\geq& r \Big(\displaystyle\frac{\partial_{t} \phi_{1} -\nabla\cdot(A\nabla\phi_{1})+q\cdot\nabla\phi_{1}}{\phi_{1}}-\mu_{1} \Big) \\
&&+ (1-r) \Big( \displaystyle\frac{\partial_{t} \phi_{2} -\nabla\cdot (A\nabla\phi_{2})+q\cdot\nabla\phi_{2}}{\phi_{2}}-\mu_{2} \Big) \\
&&\\
&\geq & r k_0(A,q,\mu_1)+(1-r)k_0(A,q,\mu_2).\\
\end{array} \]

Using the min-max characterization (\ref{minmax-dep}) of $k_{\lambda}(A,q,\mu)$, we get
\begin{equation} \label{strictconc} k_0(A,q,\mu)\geq r k_0(A,q,\mu_1)+(1-r) k_0(A,q,\mu_2).\end{equation}

This gives the first part of Proposition \ref{concavitydep}. Assume now that the equality holds. Then (\ref{strictconc0}) is an equality for all $(t,x)\in\R\times\R^N$ and thus $\nabla z_1\equiv \nabla z_2$. 
Write $z_1(t,x)=z_2(t,x)+f(t)$ for all $(t,x)\in\R\times\R^N$. Then $\phi_1(t,x)=\phi_2(t,x)e^{f(t)}$ and
\[\begin{array}{rcl}
   0&=&\partial_t\phi_1-\nabla\cdot (A\nabla\phi_1)+q\cdot \nabla\phi_1 - \mu_1\phi_1-k_0(A,q,\mu_1)\phi_1\\
&=&e^{f(t)} \Big(f'(t)\phi_2+ \partial_t\phi_2-\nabla\cdot (A\nabla\phi_2)+q\cdot \nabla\phi_2 - \mu_1\phi_2-k_0(A,q,\mu_1)\phi_2 \Big)\\
&=&e^{f(t)} \Big(f'(t)\phi_2+(\mu_2- \mu_1)\phi_2+(k_0(A,q,\mu_2)-k_0(A,q,\mu_1))\phi_2 \Big),\\
  \end{array}\]
Hence:
\[\mu_2-\mu_1\equiv -f'(t)+k_0(A,q,\mu_1)-k_0(A,q,\mu_2),\]
and the right-hand side only depends on $t$. 

\bigskip

In the other hand, if $\eta=\mu_1-\mu_2$ does not depend on $x$, set 
\[\psi(t,x)=\phi_2(t,x)\exp \Big(\int_{0}^{t}r\eta(s)ds-\frac{t}{T}\int_{0}^{T}r\eta(s)ds\Big).\]
This function is periodic in $t$ and $x$ and satisfies:
\[\begin{array}{l}\partial_{t}\psi - \nabla\cdot (A\nabla \psi)+q\cdot\nabla\psi-(\mu_2+r\eta(t))\psi=(k_{0}(A,q,\mu_2)-\frac{r}{T}\int_{0}^{T}\eta(t)dt)\psi. \\
  \end{array}
\]
The uniqueness of the eigenelements gives:
\[k_{0}(A,q,r\mu_1+(1-r)\mu_2)=k_{0}(A,q,\mu_2+r\eta)=k_{0}(A,q,\mu_2)-\frac{r}{T}\int_{0}^{T}\eta(t)dt.\]
Thus for all $r\in (0,1)$:
\[k_{0}(A,q,r\mu_1+(1-r)\mu_2)=rk_{0}(A,q,\mu_1)+(1-r)k_{0}(A,q,\mu_2).\]
$\Box$


\subsection{A general dependence result}

In \cite{Base2}, Berestycki, Hamel and Roques proved that, if $A, q$ and $\mu$ are constant, if $\eta$ does not depend on $t$ and if $\int_C\eta\geq 0$, 
then $B\mapsto k_{\lambda}(A,q,\mu+B\eta)$ is a nonincreasing function. In order to prove some of our results, we need to extend this property to general heterogeneous 
coefficients. This extension involves the principal eigenfunction $\tilde{\phi}_{\lambda}$ associated with
the adjoint problem, defined up to multiplication by a positive constant by:
\begin{equation} \left\{ \begin{array}{l}
-\partial_{t}\tilde{\phi}_{\lambda} - \nabla \cdot(A\nabla\tilde{\phi}_{\lambda}) +2\lambda e A\nabla\tilde{\phi}_{\lambda} -\nabla\cdot (q\tilde{\phi}_{\lambda})\\
 -(-\lambda\nabla\cdot (A e)+\lambda^{2} eAe-\lambda q\cdot e+\mu) \tilde{\phi}_{\lambda}=k_{\lambda}(\mu)\tilde{\phi}_{\lambda}, \\
\tilde{\phi}_{\lambda} > 0, \\
\tilde{\phi}_{\lambda} \ \hbox{is periodic}.\\
\end{array}\right. \end{equation}
We normalize this adjoint eigenfunction by $\int_{(0,T)\times C} \phi_{\lambda} \tilde{\phi}_{\lambda} = 1$. 

\begin{prop} \label{thmamplitude}

Take $\eta$ a periodic continuous function. If $\int_{(0,T)\times C} \eta\phi_{\lambda}\tilde{\phi}_{\lambda}
\geq 0$ (resp. $\int_{(0,T)\times C} \eta\phi_{\lambda}\tilde{\phi}_{\lambda}>0$), then the function $B \mapsto k_{\lambda}(A,q,\mu+B\eta)$ is
nonincreasing (resp. decreasing) over $\mathbb{R}^{+}$. 
Moreover, if $\int_{(0,T)\times C} \eta\phi_{\lambda}\tilde{\phi}_{\lambda}
= 0$ and if $A, q$ and $\mu$ do not depend on $t$, then the function $B \mapsto k_{\lambda}(A,q,\mu+B\eta)$ is
decreasing over $\mathbb{R}^{+}$ if and only if $\eta$ is not a constant with respect to $x$.
\end{prop}

\noindent {\bf Proof.}
Set $F(B)=k_{\lambda}(A,q,\mu+B\eta)$. This function is concave and analytic from the Kato-Rellich theorem. It has been proved in Theorem 3.3 of 
\cite{eigenvaluearticle} that $F'(0)=-\int_{(0,T)\times C}\mu\phi_{\lambda}\tilde{\phi}_{\lambda}dtdx$. Thus if this quantity is negative, $F$ 
is clearly decreasing over $\R^{+}$. If it is null, then $F$ is nonincreasing. As this is true for all $\lambda>0$ and as $c^*_e(A,q,\mu+B\eta)$ is a minimum, 
the opposite monotonicity properties are also true for $B\mapsto c^*_e(A,q,\mu+B\eta)$.

Moreover, if $A$, $q$ and $\mu$ do not depend on $t$, Proposition \ref{concavitydep} yields that $F$ is strictly concave if and only 
if $\eta$ is not a constant with respect to $x$. Hence, if $F'(0)=-\int_{(0,T)\times C}\eta\phi_{\lambda}\tilde{\phi}_{\lambda}dtdx=0$, then $F$ is decreasing 
if and only $\eta$ is not a constant with respect to $x$. $\Box$


\subsection{Applications of Proposition \ref{thmamplitude}}

We are now in position to prove our dependence results using Proposition \ref{thmamplitude}. 

\bigskip

\noindent {\bf Proof of Theorem \ref{influencetimevar}.}
We set $\eta=\mu-\hat{\mu}$ and we apply Proposition $\ref{thmamplitude}$, replacing $\mu$ by $\hat{\mu}$. The function $F:B\mapsto k_{\lambda}(A,q,\hat{\mu}+B\eta)$ is nonincreasing if $\int_{(0,T)\times C}\mu\phi_{\lambda}\tilde{\phi}_{\lambda}\geq 0$, where $\phi_{\lambda}$ and $\tilde{\phi}_{\lambda}$ are associated with the coefficients $(A,q,\hat{\mu})$. As $(A,q,\hat{\mu})$ do not depend on $t$, these eigenfunctions do not depend on $t$ and thus:
\[\int_{(0,T)\times C}\eta\phi_{\lambda}\tilde{\phi}_{\lambda}=\int_{C}(\int_{0}^{T}(\mu(t,x)-\hat{\mu}(x))dt)\phi_{\lambda}(x)\tilde{\phi}_{\lambda}(x)dx=0\]
since $\hat{\mu}(x)=\frac{1}{T}\int_{0}^{T}\mu(t,x)dt$. Thus:
\[F(1)=k_{\lambda}(A,q,\hat{\mu}+\eta)=k_{\lambda}(A,q,\mu)\leq F(0)= k_{\lambda}(A,q,\hat{\mu}).\]
As this is true for all $\lambda\in\R^{+}$, this gives:
\[c^*_e(A,q,\mu)\geq c^*_e(A,q,\hat{\mu}).\]

If the equality holds, considering some $\lambda\in\R^N$ such that $\lambda\cdot e<0$ and $c^*_e(A,q,\mu)=\frac{k_{\lambda}(A,q,\mu)}{\lambda\cdot e}$, one gets 
\[\frac{k_{\lambda}(A,q,\mu)}{\lambda\cdot e}=c^*_e(A,q,\mu)=c^*_e(A,q,\hat{\mu})\leq \frac{k_{\lambda }(A,q,\hat{\mu})}{\lambda\cdot e}.\]
Thus $k_{\lambda}(A,q,\hat{\mu}+\eta)=k_{\lambda}(A,q,\hat{\mu})$. Proposition \ref{thmamplitude} then gives that $\eta$ does not depend on $x$. 
Thus $\mu(t,x)=\hat{\mu}(x)+\eta(t)$. $\Box$

\bigskip

\noindent {\bf Proof of Proposition \ref{dep-x}.} First of all, if $A$, $q$ and $\mu$ do not depend on $x$, using Proposition 3.1 
of \cite{eigenvaluearticle}, we have:
\[k_{\lambda}(A,q,\mu)=-\frac{1}{T}\int_{0}^{T}(\lambda A \lambda-\lambda\cdot q+\mu).\]
In order to compute $\min_{\lambda\cdot e<0}\frac{1}{\lambda\cdot e}\frac{1}{T}\int_{0}^{T}(\lambda A \lambda-\lambda\cdot q+\mu)$, let write
$\lambda=\alpha \xi$, with $\xi \in\mathbb{S}^{N-1}$ and $\alpha<0$. We compute:
$$\begin{array}{rcl}
\min_{\lambda\cdot e<0}\frac{1}{\lambda\cdot e}\frac{1}{T}\int_{0}^{T}(\lambda A \lambda-\lambda\cdot q+\mu)&=&\min_{\xi\cdot e>0}\min_{\alpha<0}\alpha\frac{1}{T}\int_{0}^{T}\xi A \xi-\frac{1}{T}\int_{0}^{T}\xi\cdot q+\frac{1}{\alpha}\frac{1}{T}\int_{0}^{T}\mu)\\
&=& \min_{e\cdot \xi>0} \frac{2}{T}\sqrt{\int_{0}^{T}\xi A \xi\int_{0}^{T}\mu}-\frac{1}{T}\int_{0}^{T}q\cdot \xi.\\
\end{array}$$
This gives the equality in (\ref{dep-x-ineq}). 

Assume now that $\mu$ depends on $x$. Set $\eta=\mu-\overline{\mu}$. The same arguments as in the proof of Theorem \ref{influencetimevar} give
$c^*_e(A,q,\mu)\geq c^*_e(A,q,\overline{\mu})$ and the equality holds if and only if $\eta$ does not depend on $x$. 
In this case, $\mu=\overline{\mu}+\eta$ does not depend on $x$. $\Box$

\bigskip

\noindent {\bf Proof of Proposition \ref{dep-amp-classical}.} 
1. This is an immediate consequence of Proposition \ref{thmamplitude} since, when $A=I_N$, $q\equiv 0$ and $\mu$ is a constant, one has 
$\phi_\lambda\equiv \tilde{\phi}_\lambda\equiv 1$.

2. Set $\hat{\eta}(x)=\frac{1}{|C|}\int_C \eta(t,x)dx$. We know from Theorem \ref{influencetimevar} that
$$k_{\lambda}(A,q,\mu+B\eta)\leq k_{\lambda}(A,q,\mu+B\hat{\eta}).$$
Moreover, Berestycki, Hamel and Roques have proved in \cite{Base2} that, as $\max_{x\in\R^N} \hat{\eta}>0$, the right-hand side goes to $-\infty$ as $B\to+\infty$. 
Thus the left-hand side converges to $-\infty$ as $B\to +\infty$. As it is a concave function of $B$, it is decreasing over $[B_0,\infty)$, with $B_0$ large enough. 
Hence $B\mapsto c^*_e(A,q,\mu+B\eta)$ is increasing. 
$\Box$


\section{Proof of the monotonicity with respect to the diffusion term}\label{section-proofsA}

\noindent {\bf Proof of Theorem \ref{mon-diffusion}.} fix $e \in\mathbb{S}^{N-1}$ and $\lambda>0$. It has been proved by the author in \cite{Steiner} that 
$$k_{\lambda e}(\kappa A,0,\mu)=\min_{\alpha\in \mathcal{A}} \Big(\int_C \kappa\nabla\alpha A(x)\nabla \alpha -
\int_C \mu(x)\alpha^2 -\lambda^2\kappa |C| D_e (\alpha^2 A) \Big),$$
where 
\[\mathcal{A}=\{ \alpha\in \mathcal{C}^1_{per}(\R^N), \ \alpha>0, \ \int_C \alpha^2=1\}.\]
and $D_e(A)$ is the effective diffusivity of a matrix field $A$ in direction $e$, that is,
\begin{equation} \label{effectivediffusivity}  D_e(A)=\min_{\chi \in\mathcal{C}^1_{per}(\R^N)} \frac{1}{|C|}\int_{C} (e+\nabla\chi)A(x)(e+\nabla\chi). \end{equation}
This formula yields that $\kappa\mapsto k_{\lambda e}(\kappa A, 0, \mu) $ is a concave function. 
In the other hand, it has been proved by Pinsky \cite{Pinsky} that $\lambda\mapsto k_{\lambda e}(A,0,\mu)$ is strictly concave and we know (see \cite{Base1, eigenvaluearticle} for 
example) that this function reaches iots maximum when $\lambda=0$, thus
$$k_{\lambda}(\kappa A, 0, \mu)< k_0(\kappa A, 0, \mu) \hbox{ for all } \lambda\in\R^N\backslash \{0\}.$$
Moreover, $k_0(\kappa A, 0,\mu)\leq -\max_{\R^N}\mu$. We now prove that $k_{\lambda e}(\kappa A, 0, \mu)\to -\max_{\R^N}\mu$ as $\kappa \to 0$ in order to conclude. 
To do so, we use Lemma \ref{compjlambda} (see below) to get:
$$k_\lambda(\kappa A, 0, \mu)\geq k_0(\kappa A, 0,\lambda A \lambda+\mu)\geq k_0(\kappa A, 0,\mu)-\gamma \kappa |\lambda|^2,$$
where $\gamma$ is the ellipticity constant given by (\ref{ellipticity-dep}). This leads to
$$ \liminf_{\kappa\to 0} k_{\lambda e}(\kappa A, 0, \mu)\geq \liminf_{\kappa\to 0} k_0(\kappa A, 0,\mu)=-\max_{\R^N}\mu,$$
where the convergence in the right-hand side has been proved by the author in \cite{eigenvaluearticle}. 
This gives $k_{\lambda e}(\kappa A, 0, \mu)\to -\max_{\R^N}\mu$ as $\kappa \to 0$

Hence, for all $\kappa >0$, $k_{\lambda e}(\kappa A, 0, \mu)\leq \lim_{\kappa'\to 0} k_{\lambda e} (\kappa' A , 0, \mu)$. 
As $\kappa\mapsto k_{\lambda e}(\kappa A, 0, \mu) $ is concave, it is then a decreasing function. This gives the conclusion. $\Box$


\section{Proof of the dependence results with respect to the drift term}\label{section-proofsq}

{\bf Proof of Proposition \ref{shearflowprop}}

This proposition relies on the following observation, which holds for general diffusion matrix $A$, drift term $q$ and reaction term $\mu$ that can be written 
$A=a(t,y)I_{N}$, $q(t,x)=(q_{1}(t,y),0,...,0)$ and $\mu=\mu(t,y)$. In this case, for all direction of propagation $e=(e_{1},\tilde{e})\in\mathbb{S}^{N-1}$, define
$\tilde{k}_\lambda(a,q_1,\mu)$ the periodic principal eigenvalue defined by the existence of a function 
$\varphi_\lambda\in\mathcal{C}_{per}^{1,2}(\R\times\R^{N-1})$ that solves
\begin{equation} \label{shearflow}\left\{ \begin{array}{l}
                          \partial_{t}\varphi_{\lambda}-\nabla \cdot (a(t,y)\nabla\varphi_{\lambda})-2\lambda a(t,y)\tilde{e}\cdot\nabla\varphi_{\lambda}\\
-(\lambda\nabla\cdot (a(t,y)\tilde{e})+\lambda^{2}a(t,y)-\lambda q_{1}(t,y)e_{1}+\mu(t,y))\varphi_{\lambda}= \tilde{k}_{\lambda}(a,q_1,\mu)\varphi_{\lambda} \ \hbox{in} \ \R\times\R^{N-1},\\
\varphi_{\lambda}>0 \ \hbox{in} \ \R\times\R^{N-1},\\
\varphi_{\lambda} \ \hbox{is periodic in $t$ and $y$}.\\
                         \end{array}\right.\end{equation}
Setting $\psi(t,x_1,y)=\varphi_\lambda(t,y)$, this function satisfies the eigenvalue problem (\ref{eigentw2}) associated with $L_\lambda$. The uniqueness of 
the periodic principal eigenvalue yields that $k_\lambda(aI_N,q,\mu)=\tilde{k}_\lambda(a,q_1,\mu)$. 
Thus if $a\equiv 1$ and $\mu$ is a positive constant, we immediatly get from the proof of Proposition \ref{dep-amp-classical} that
\[B\mapsto k_{\lambda}(I_N,Bq,\mu_0) \hbox{ is decreasing},\]
which conludes the proof. $\Box$

\bigskip

The proof of Theorem \ref{potdrift} uses a lemma of independent interest that we state separately:

\begin{lem} \label{compjlambda}
For all coefficients $(A,q,\mu)$ that do not depend on $t$ and $\lambda\in\R^N$, one has
\[k_\lambda(A,q,\mu) \geq k_0(A,0,\frac{\nabla\cdot q}{2}+\lambda A \lambda-\lambda \cdot q+\mu).\]
\end{lem}

\noindent {\bf Proof.} We know that $L_\lambda\phi_\lambda=k_\lambda (A,q,\mu)\phi_\lambda$, where $\phi_\lambda$ does not depend on $t$. Multiplying this equation by $\phi_\lambda$ and integrating over $C$, this gives:
\[\begin{array}{rcl}  k_\lambda(A,q,\mu) \int_{C}\phi_\lambda^2&=&\int_{C}\nabla\phi_\lambda A \nabla\phi_\lambda + \frac{1}{2}\int_{C} (q-2A\lambda)\nabla (\phi_\lambda^2) \\
&&- \int_{ C} (\lambda A\lambda +\nabla\cdot (A\lambda) -q\cdot\lambda +\mu)\phi_\lambda^2\\
&&\\
&=&\int_{ C}\nabla\phi_\lambda A \nabla\phi_\lambda-\int_{C} (\lambda A\lambda +\frac{\nabla\cdot q}{2}-q\cdot\lambda +\mu)\phi_\lambda^2. \\ \end{array} \]
But we know from the Rayleigh characterization that if
\[X= \{\phi\in\mathcal{C}^{2}_{per}(\R^N), \ \phi>0, \ \int_{C}\phi^2=1 \},\]
then:
\[\begin{array}{l}  k_0(A,0,\frac{\nabla\cdot q}{2}+\lambda A \lambda-\lambda \cdot q+\mu)\\
= \min_{\phi\in X}\int_{C}\nabla\phi A \nabla\phi-\int_{(0,T)\times C} (\lambda A\lambda +\frac{\nabla\cdot q}{2}-q\cdot\lambda +\mu)\phi^2.\\ \end{array} \]
This gives the conclusion. $\Box$

\bigskip

\noindent {\bf Proof of Theorem \ref{potdrift}.} As $\nabla Q$ is periodic and $\int_C q =0$, it has been proved in \cite{Rossi} that $Q$ is periodic. 
Set $\phi_{\lambda}$ a positive eigenfunction associated with
$k_{\lambda }$ and $\psi_{\lambda}(x)=\phi_{\lambda}(x)e^{-Q(x)/2}$.
This new function satisfies:
\begin{equation} \left\{ \begin{array}{l}
L_{\lambda}\psi_{\lambda} -(\frac{1}{2}\Delta Q-\frac{1}{4}|\nabla Q|^2) \psi_{\lambda} = k_{\lambda}\psi_{\lambda}, \\
\psi_{\lambda} > 0, \\
\psi_{\lambda} \ \hbox{is periodic}.\\
\end{array} \right. \end{equation}
This yields that:
\[k_{\lambda}(I_N,\nabla Q,\mu)=k_{\lambda}(I_N,0,\frac{1}{2}\Delta Q-\frac{1}{4}|\nabla Q|^2+\mu),\]
which gives 
\[c^*_e(I_N,\nabla Q,\mu)=c^*_e(I_N,0,\mu-\frac{1}{4}|\nabla Q|^ 2-\Delta Q)).\]

\smallskip

We are now in position to prove the results. Using Lemma \ref{compjlambda}, we get
\begin{equation} \begin{array}{rcl}
 \displaystyle\frac{1}{\lambda B^{2}}k_{\lambda B}(I_N,B\nabla Q,\mu_0)&=&\displaystyle\frac{1}{\lambda B^{2}}k_{\lambda B}(I_N,0,\displaystyle\frac{B}{2}\Delta Q-\displaystyle\frac{B^{2}}{4}|\nabla Q|^2+\mu_0)\\
\\
&\geq &\displaystyle\frac{ 1}{\lambda B^{2}}k_{0}(I_N,0,\lambda^{2}B^{2}-\displaystyle\frac{B}{2}\Delta Q-\displaystyle\frac{B^{2}}{4}|\nabla Q|^2+\mu_0)\\
\\
&\geq & \displaystyle\frac{1}{\lambda B^{2}} k_{0}(I_N,0,\lambda^{2}B^{2}-\displaystyle\frac{B^{2}}{4}|\nabla Q|^2)-\|\displaystyle\frac{B}{2}\Delta Q+\mu_0\|_{\infty}\displaystyle\frac{1}{B^{2}} \\
\\
&\geq& -\lambda+\displaystyle\frac{1}{\lambda B^{2}}k_{0}(I_N,0,-\displaystyle\frac{B^{2}}{4}|\nabla Q|^2)+O(1/B)\geq -\lambda+O(1/B).\\ \end{array}\end{equation}
Thus, one gets
\begin{equation} \begin{array}{rcl}\liminf_{B\rightarrow+\infty}\displaystyle  \displaystyle\frac{k_{\lambda B}(I_N,B\nabla Q,\mu_0)}{\lambda B^{2}}&\geq& -\lambda. \\
\end{array}\end{equation}
In the other hand, we know that for all $\lambda>0$, one has
\[\frac{1}{B}c^*_e(I_N,B\nabla Q,\mu_0)\leq \frac{-k_{\lambda B}(I_N,B\nabla Q,\mu_0)}{\lambda B^{2}}.\]
Letting $B\to+\infty$, this gives
$$\limsup_{B\to +\infty} \frac{1}{B}c^*_e(I_N,B\nabla Q,\mu_0) \leq \lambda.$$
As this is true for all $\lambda>0$, one has $\limsup_{B\to +\infty} \frac{1}{B}c^*_e(I_N,B\nabla Q,\mu_0) \leq 0.$

Next, for all $\lambda$, it has been proved in \cite{Base1} that
$$\begin{array}{rcl}
k_\lambda(I_N,B\nabla Q,\mu_0)&=&k_{\lambda}(I_N,0,\displaystyle\frac{B}{2}\Delta Q-\displaystyle\frac{B^{2}}{4}|\nabla Q|^2+\mu_0)\\
&&\\
&\leq& k_0(I_N,0,\displaystyle\frac{B}{2}\Delta Q-\displaystyle\frac{B^{2}}{4}|\nabla Q|^2+\mu_0)\\
&&\\
&\leq&k_0(I_N,B\nabla Q,\mu_0)=-\mu_0<0.\\
\end{array}$$
Thus $\limsup_{B\to +\infty} \frac{1}{B}c^*_e(I_N,B\nabla Q,\mu_0) \geq 0$, which proves 2. of Theorem \ref{potdrift}.

\smallskip

If $\mu_0$ is a positive constant and $Q$ does not depend on $t$, using Lemma \ref{compjlambda}, we get:
\begin{equation} \begin{array}{rcl}
k_{\lambda}(I_N,B\nabla Q,\mu_0)&=& k_{0}(I_N,-2\lambda e,|\lambda|^{2}-\displaystyle\frac{B^{2}}{4}|\nabla Q|^2+\displaystyle\frac{B}{2}\Delta Q+\mu_0) \\
&&\\
&\geq & k_{0}(I_N,0,|\lambda|^{2}-\displaystyle\frac{B^{2}}{4}|\nabla Q|^2+\displaystyle\frac{B}{2}\Delta Q+\mu_0)\\
&&\\
&&= k_0(I_N,B\nabla Q, |\lambda|^2 +\mu_0)=-|\lambda|^2 -\mu_0.\\
\end{array}\end{equation}
This proves 2. of Theorem \ref{potdrift}.
$\Box$


\begin{thebibliography}{10}


\bibitem{Aronson}
D.G. Aronson and H.F. Weinberger.
\newblock Multidimensional nonlinear diffusions arising in population genetics.
\newblock {\em Adv. Math.}, 30:33--76, 1978.

\bibitem{Audoly}
B.~Audoly, H.~Berestycki, and Y.~Pomeau.
\newblock R\'eaction diffusion en \'ecoulement stationnaire rapide.
\newblock {\em C. R. Acad. Sci. Paris}, 328:255--262, 2000.

\bibitem{BHNa}
H.~Berestycki, F.~Hamel, and G.~Nadin.
\newblock Asymptotic spreading in heterogeneous diffusive excitable media.
\newblock {\em J. Func. Anal.}, 255(9):2146--2189, 2008.

\bibitem{Frontexcitable}
H. Berestycki and F. Hamel.
\newblock Front propagation in periodic excitable media.
\newblock {\em Comm. Pure Appl. Math.}, 55:949--1032, 2002.

\bibitem{Largedrift}
H.~Berestycki, F.~Hamel, and N.~Nadirashvili.
\newblock Elliptic eigenvalue problems with large drift and applications to
  nonlinear propagation phenomena.
\newblock {\em Comm. Math. Phys.}, 253:451--480, 2005.

\bibitem{SpeedKPP}
H.~Berestycki, F.~Hamel, and N.~Nadirashvili.
\newblock The speed of propagation for kpp type problems. i - periodic
  framework.
\newblock {\em J. Europ. Math. Soc.}, 7:173--213, 2005.

\bibitem{Base1}
H.~Berestycki, F.~Hamel, and L.Roques.
\newblock Analysis of the periodically fragmented environment model.
              {I}. {S}pecies persistence.
\newblock {\em J. Math. Biol.}, 51:75--113, 2005.

\bibitem{Base2}
H.~Berestycki, F.~Hamel, and L.Roques.
\newblock Analysis of the periodically fragmented environment model : {II} -
  {B}iological invasions and pulsating travelling fronts.
\newblock {\em J. Math. Pures Appl.}, 84:1101--1146, 2005.

\bibitem{Rossi}
H. Berestycki, F. Hamel and L. Rossi.
\newblock Liouville-type results for semilinear elliptic equations in
              unbounded domains.
\newblock {\em Ann. Mat. Pura Appl.}, 186 (4):469--507, 2007.

\bibitem{ElSmaily}
M.~ElSmaily.
\newblock Pulsating travelling fronts: Asymptotics and homogenization regimes.
\newblock {\em preprint}, 2007.

\bibitem{Fisher}
R.~A. Fisher.
\newblock The advance of advantageous genes.
\newblock {\em Ann. Eugenics}, 7:335--369, 1937.

\bibitem{Freidlin2}
M.~Freidlin.
\newblock On wave front propagation in periodic media.
\newblock {\em In: Stochastic analysis and applications, ed. M. Pinsky,
  Adavances in Probability and related topics}, 7:147--166, 1984.

\bibitem{Gartner}
M.~Freidlin and J.~Gartner.
\newblock On the propagation of concentration waves in periodic and random
  media.
\newblock {\em Sov. Math. Dokl.}, 20:1282--1286, 1979

\bibitem{HamelRoqueslimgamma}
F. Hamel, L. Roques, and J. Fayard.
\newblock Spreading speeds in slowly oscillating environments. 
\newblock {\em preprint}, 2009.

\bibitem{limgamma}
F. Hamel, G. Nadin, and L. Roques.
\newblock ???. 
\newblock {\em preprint}, 2009.

\bibitem{Heinze2005}
S.~Heinze.
\newblock Large convection limits for kpp fronts.
\newblock {\em Max Planck Institute for Mathematics Preprint Nr.}, 2005.

\bibitem{Heinze}
S.~Heinze, G.~Papanicolaou, and A.~Stevens.
\newblock Variational principles for propagation speeds in inhomogeneous media.
\newblock {\em SIAM J. Appl. Math.}, 62(1):129--148, 2001.

\bibitem{Hutson}
V.~Hutson, K.~Michaikow, and P.~Polacik.
\newblock The evolution of dispersal rates in a heterogeneous time-periodic
  environment.
\newblock {\em J. Math. Biol.}, 43:501--533, 2001.

\bibitem{Hutson2}
V.~Hutson, W.~Shen, and G.~T. Vickers.
\newblock Estimates for the principal spectrum point for certain time-dependent
  parabolic operators.
\newblock {\em Proc. A.M.S.}, 129:1669--1679, 2000.

\bibitem{Kiselevadvection}
A.~Kiselev and L.~Ryzhik.
\newblock Enhancement of the traveling front speeds in reaction-diffusion
  equations with advection.
\newblock {\em Ann. Inst. H. Poincar\'e Anal. Non Lin\'eaire}, 18:309--358,
  2001.

\bibitem{KPP}
A.N. Kolmogorov, I.G. Petrovsky, and N.S. Piskunov.
\newblock Etude de l \'equation de la diffusion avec croissance de la
  quantit\'e de mati\`ere et son application \`a un probl\`eme biologique.
\newblock {\em Bulletin Universit\'e d'Etat \`a Moscou (Bjul. Moskowskogo Gos.
  Univ.)}, pages 1--26, 1937.

\bibitem{eigenvaluearticle}
G.~Nadin.
\newblock The principal eigenvalue of a space-time periodic parabolic operator.
\newblock {\em Ann. Mat. Pura Appl.}, 4:269--295, 2009.

\bibitem{twperiodicarticle}
G.~Nadin.
\newblock Traveling fronts in space-time periodic media.
\newblock {\em to appear in J. Math. Pures Appl.}, 2009.

\bibitem{Steiner}
G.~Nadin. 
\newblock The effect of the Schwarz rearrangement on the periodic principal eigenvalue of a nonsymmetric operator.
\newblock {\em to appear in SIAM J. Math. Anal.}, 2009.

\bibitem{NolenRuddXin}
J.~Nolen, M.~Rudd, and J.~Xin.
\newblock {Existence of {KPP} fronts in spatially-temporally periodic advection
  and variational principle for propagation speeds}, Dynamics of
  PDE.
\newblock 2(1):1--24, 2005.

\bibitem{Nolenshear}
J.~Nolen and J.~Xin.
\newblock Existence of {KPP} type fronts in space-time periodic shear flows and a
  study of minimal speeds based on variational principle.
\newblock {\em Disc. and Cont. Dyn. Syst.}, 13(5):1217--1234, 2005.

\bibitem{NolenXin1d}
J.~Nolen and J.~Xin.
\newblock Kpp fronts in 1d random drift.
\newblock {\em Discrete and Continuous Dynamical Systems B, accepted for
  publication}, 2008.

\bibitem{PapanicolaouXin}
G.~Papanicolaou and X.~Xin.
\newblock Mathematical biology.
\newblock {\em J. Stat. Phys.}, 63:915--932, 1991.

\bibitem{Pinsky}
R. G. Pinsky.
\newblock Second order elliptic operators with periodic coefficients: criticality theory, perturbations, and positive harmonic functions.
\newblock {\em J. Funct. Anal.}, 129:80--107, 1995.

\bibitem{RyzhikZlatos}
L.~Ryzhik and A.~Zlatos.
\newblock {KPP} pulsating front speed-up by flows.
\newblock {\em Commun. Math. Sci.}, 5:575--593, 2007.

\bibitem{Weinberger}
H.~Weinberger.
\newblock On spreading speed and travelling waves for growth and migration
  models in a periodic habitat.
\newblock {\em J. Math. Biol.}, 45:511--548, 2002.

\bibitem{Xin}
J. Xin,
\newblock Existence of planar flame fronts in convective-diffusive periodic media.
\newblock {\em Arch. Ration. Mech. Anal.}, 121:205--233, 1992.

\bibitem{Zlatosc/A}
A.~Zlatos.
\newblock Sharp asymptotics for {kpp} pulsating front speed-up and diffusion
  enhancement by flows.
\newblock {\em to appear in Arch. Ration. Mech. Anal.}, 2009.

\bibitem{Zlatoscex}
A.~Zlatos.
\newblock Reaction-diffusion front speed enhancement by flows.
\newblock {\em preprint}, 2009.


\end{thebibliography}
\end{document}